\documentclass[a4paper,12pt]{article}
\usepackage[utf8]{inputenc}
\usepackage[english]{babel}
\usepackage{amsmath}
\usepackage{amssymb}
\usepackage{graphicx}
\usepackage{caption}
\usepackage{subcaption}

\graphicspath{{pictures/}}
\DeclareGraphicsExtensions{.pdf,.png,.jpg}

\title{Subresonant solutions of the linear oscillator equation}
\author{P.Y.Astafyeva, O.M.Kiselev}
\begin{document}
	\maketitle
	\begin{abstract}
		The behavior of a linear oscillator under the action of an external almost periodic force is investigated. The constructed solutions grow more slowly than the resonant ones. The dependence of the amplitude of growing solutions on the parameters of an almost periodic perturbation is calculated.
	\end{abstract}
	
\section{Introduction}
\label{secIntroduection}

	In this paper, the behavior of a linear oscillator under the action of an external force is considered:
	\begin{equation}
		\frac{d^2}{dt^2}u+u=f(t),
		\label{eqLinOsc1}
	\end{equation}
where $f(t)$ as an almost periodic function:
\begin{equation}
f=\sum_{n=1}^{\infty}\frac{1}{n^k}\cos((1-\frac{1}{n^p})t).
\label{formulaForF}
\end{equation}

Two types of solutions are well studied among the solutions of the linear oscillator equation with an oscillating external force. The first type of solutions arises if the frequency of oscillations of the external force does not coincide with the natural oscillator. An example of such solutions is an equation with a periodic right-hand side of the following form:
	$$f=a \cos(\phi +k t ).$$
In this case, the general solution is:
$$u=C_{1}\cos(t+\alpha)  +\frac{a \cos( \phi +k t)}{1-k^2},$$	
 here $k\not=\pm1$,  $ C_{1}, \alpha$ -- solution parameters. If $k\not\in\mathbb{N}$, when solutions are run-out. If $k\in\mathbb{N}$, when solution-- periodic function with period $T=2\pi k $.
 
 If $k=1$, then the amplitude of oscillations of the solution grows linearly. The general solution is:
 $$u=C_{1}\cos(t+\alpha)+\frac{a t \sin( t +\phi)}{2}$$ 
This solution is called resonant.

Exact coincidence of the frequencies of the external force and the oscillator is a rare phenomenon in a general situation. Exact frequency matching cannot be achieved with frequency tuning. For example, This doesn't done in local oscillators - automatic frequency control systems. As a rule, the result is a frequency close to resonant, however, this frequency is, strictly speaking, unstable.

It is important to consider the behavior of the oscillator
under the action of an external force with a close to resonance
frequency and, moreover, with a set of close frequencies in such conditions. External forces with almost periodic behavior are typical here. General approaches to the study of almost-periodic functions are presented, for example, in \cite{Levitan1953}, \cite{LevitanZhikov1977}. 

The paper shows the effect on a linear oscillator of an almost periodic force having a set of frequencies, among which there are close to resonant ones, but their amplitude in the general set of frequencies of an almost periodic perturbation tends to
zero. The amplitude of the forced oscillations of the linear oscillator under the action of such almost-periodic perturbations increases more slowly than in the case of resonance. Such solutions are called subresonant in the article.

Close to resonant solutions lead to small denominators. Small denominators appear as defects in the asymptotic formalism and, as a rule, can be eliminated by modulating the parameters of asymptotic expansions in perturbation theory \cite{Arnold1963}. In the proposed paper, small denominators arise essentially due to the unavoidable proximity of the set of frequencies of the external force to the resonant frequency. 

In section \ref{secFormalSolution} formulas for solving the equation are shown (\ref{eqLinOsc1}) with the right side (\ref{formulaForF}). In section \ref{secSymptoticsOfSeries} asymptotics by $t$ are obtained for the series appearing in the formula for the solution. In section \ref{secAsymptoticsOfSubresonantSolution} the asymptotic behavior of the subresonant solution is presented, and the main result is formulated.

\section{Formal solution of an equation with almost periodic perturbation}
\label{secFormalSolution}

We will seek a solution to the Cauchy problem for an equation (\ref{eqLinOsc1}) with zero initial conditions:
\begin{eqnarray*}
	\left\{
	\begin{array}{c}	
		u_n|_{t=0}=0
		\\
		u'_n|_{t=0}=0
	\end{array}
	\right.
\end{eqnarray*}
The solution to this problem can be represented as:
 \begin{eqnarray} u&=&\frac{1}{2}\sin(t)\sum_{n=1}^{\infty}(-\frac{n^p\sin(\frac{t}{n^p}-2t)}{2n^{p+k}-n^k}+\frac{2n^{2p}\sin(\frac{t}{n^p})}{2n^{p+k}-n^k}-\frac{n^p\sin(\frac{t}{n^p})}{2n^{p+k}-n^k})+\nonumber\\&+&\frac{1}{2}\cos(t)\sum_{n=1}^{\infty}(\frac{n^p\cos(\frac{t}{n^p}-2t)}{2n^{p+k}-n^k}+\frac{2n^{2p}(\cos(\frac{t}{n^p})-1)}{2n^{p+k}-n^k}-\frac{n^p\cos(\frac{t}{n^p})}{2n^{p+k}-n^k}).
 	\label{resh}
 \end{eqnarray}
 The amplitude of oscillation is determined coefficients of $\sin(t) $ and $\cos(t) $.
 We select the uniformly bounded functions in $ t $ from the right-hand side of the formula. The first and third terms  in the series for $\sin(t)$ and $\cos(t)$  are denoted respectively as follows:
 \begin{eqnarray*}
 	R_{1}(t)=\frac{1}{2}\sin(t)\sum_{n=1}^{\infty}(-\frac{n^p\sin(\frac{t}{n^p}-2t)}{2n^{p+k}-n^k}-\frac{n^p\sin(\frac{t}{n^p})}{2n^{p+k}-n^k});
 	\\ R_{2}(t)=\frac{1}{2}\cos(t)\sum_{n=1}^{\infty}(\frac{n^p\cos(\frac{t}{n^p}-2t)}{2n^{p+k}-n^k}-\frac{n^p\cos(\frac{t}{n^p})}{2n^{p+k}-n^k}).
 \end{eqnarray*}
Using the integral convergence criterion, we can show, that the series $R_{1}(t)$ and $R_{2}(t)$  are uniformly bounded in $ t $, for $k>1$.

It remains to investigate the behavior with respect to $ t $ of two series:
\begin{equation}
R_3=\sum_{n=1}^{\infty}\frac{2n^{2p}\sin(\frac{t}{n^p})}{2n^{p+k}-n^k},
\quad
R_4=\sum_{n=1}^{\infty}\frac{2n^{2p}(\cos(\frac{t}{2n^p})-1)}{2n^{p+k}-n^k}.
\label{formulasForR3R4}
\end{equation}

Using:
$$\cos\left(\frac{t}{n^p}\right)-1=-2\sin^2\left(\frac{t}{2n^p}\right)$$
 the row $ R_4 $ can be represented as:
\begin{equation}
R_4=-2\sum_{n=1}^{\infty}\frac{2n^{2p}\sin^{2}(\frac{t}{2n^p})}{2n^{p+k}-n^k}
\label{fomulaForR4InSineForm}
\end{equation}

The coefficient in the common term of the $ R_3 $ and $ R_4 $ series can be conveniently written as:
\[
K_{n}=\frac{{{n}^{p-k}}}{1-\frac{1}{{{2n}^{p}}}}=n^{p-k}\sum_{m=0}^{\infty}(\frac{1}{2n^p})^m.
\]
Then:
\begin{eqnarray*}
R_3
&=&
\sum_{n=1}^{\infty}n^{p-k}\sin(\frac{t}{n^p})\sum_{m=0}^{\infty}(\frac{1}{2n^p})^m=
\\
&\,&
\sum_{n=1}^{\infty}n^{p-k}\sin(\frac{t}{n^p})+\sum_{n=1}^{\infty}n^{p-k}\sin(\frac{t}{n^p})\sum_{m=1}^{\infty}(\frac{1}{2n^p})^m.
\end{eqnarray*}

The second sum can be rewritten as:
\[
\left|\sum_{n=1}^{\infty}n^{p-k}\sin(\frac{t}{n^p})\sum_{m=1}^{\infty}(\frac{1}{2n^p})^m\right|\le\frac{1}{2}\sum_{n=1}^{\infty}n^{-k}\sum_{m=0}^{\infty}\frac{1}{2^{m}n^{pm}}
\]
The series converge on the right-hand side of the inequality. Similarly, you can convert the series to $ R_4 $.

It remains to investigate the behavior of the series:

\[\sigma_{s}=\sum_{n=1}^{\infty }{\left. {{n}^{p-k}} \sin{\left( \frac{t}{{{n}^{p}}}\right) }\right.}\]
\[\sigma_{c}=\sum_{n=1}^{\infty }{\left. {{n}^{p-k}}\,  \sin^{2}{\left( \frac{t}{{{2n}^{p}}}\right) } \right.}\]
Consider the series $\sigma_{s}$. Let us show that it converges absolutely for any value of $ t $. To do this, we expand the sine function in a Maclaurin series, as a result we get:
$$ \sigma_{s}= \sum_{n=1}^{\infty }{\left. {{n}^{p-k}}\, \sum_{j=1}^{\infty }{(-1)^{j-1}\left. \frac{{{t}^{2 j-1}}}{\left( 2 j-1\right) ! {{n}^{\left( 2 j-1\right)  p}}}\right.}\right.} $$
In this formula, we will change the order of summation.

$$ \sigma_{s}=\sum_{j=1}^{\infty }(-1)^{j-1} \frac{t^{2 j-1}}{( 2 j-1) ! } \sum_{n=1}^{\infty } \frac{n^{p-k}}{n^{( 2 j-1)  p}}=\sum_{j=1}^{\infty }(-1)^{j-1} \frac{t^{2 j-1}}{( 2 j-1) ! } \sum_{n=1}^{\infty } \frac{1}{n^{ 2 jp+k}} $$

Then: 
$$\sigma_{s}=\sum_{j=1}^{\infty }(-1)^{j-1} \frac{t^{2 j-1}}{( 2 j-1)!}\zeta(2jp+k)$$
The function $\zeta(2jp + k)$ is uniformly bounded in $j\in\mathbb{N}$, for $ p> 0 $, $ k> 1 $. Then $\sigma_{s} $ can be majorized by series:
$$ |\sigma_{s}|\le C\sum_{j=1}^{\infty } \frac{t^{2 j-1}}{( 2 j-1)!}, $$
   where $C>\zeta(k)$.  The series on the right side converges for any $ t $.
   The same estimates are valid for the series $\sigma_{c} $.

\section{Asymptotics of the series $\sigma_{s}$ as $t\rightarrow \infty$}
\label{secSymptoticsOfSeries}

\begin{figure}[t]
	 \includegraphics[scale=0.5]{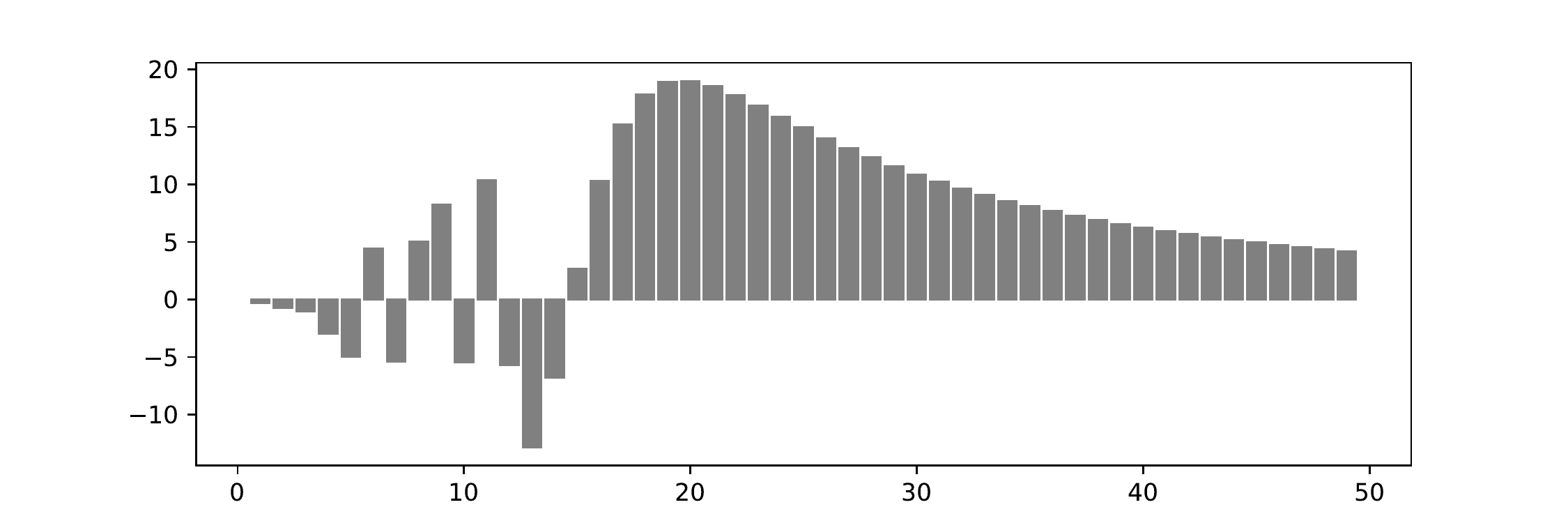} 
\caption{Typical form of a function under the sum sign. The function sign is variable up to $n=\big[(2t/\pi)^{1/p}\big]$, where the square brackets represent the integer part of the number. Here $t=10000$, $k=2$, $p=3$.}
\label{fig1}
\end{figure}

We denote $N=\big[(2t/\pi)^{1/p}\big]$, where square brackets denote the integer part of the number. We represent the series as the sum of a segment with an oscillating part and a series of fixed sign:
\begin{equation}
\sigma_{s}=\sum_{n=1}^{N }{\left. {{n}^{p-k}} \sin{\left( \frac{t}{{{n}^{p}}}\right) }\right.}+\sum_{n=N+1}^{\infty}{\left. {{n}^{p-k}} \sin{\left( \frac{t}{{{n}^{p}}}\right) }\right.}
\label{formulasigmaS}
\end{equation}

For the series on the right side of the formula (\ref{formulasigmaS}), write the majorant and minorant:
\begin{equation}
	\int_{N+1}^{\infty}{{n}^{p-k}} \sin\left( \frac{t}{n^{p}}\right) dn
	\le
	\sum_{n=N+1}^{\infty} {n^{p-k}} \sin\left( \frac{t}{{n^{p}}}\right)
	\le
	\int_{N}^{\infty}{{n}^{p-k}} \sin\left( \frac{t}{n^{p}}\right) dn 
\end{equation}
Let us estimate the asymptotics of the integral over $ t $. For this, it is convenient to make the replacement in the integrals.
Then we get:
\begin{equation}
\int_{N}^{\infty}{{n}^{p-k}} \sin\left( \frac{t}{n^{p}}\right) dn=\frac{{t}^{1-\frac{k-1}{p}}}{p} \int_{0}^{\frac{\ensuremath{\pi} }{2}}{\left. {{\tau }^{\frac{k-1}{p}-2}} \sin{\left( \tau \right) }\right.}d\tau 
\label{integral}
\end{equation}
The integral on the right side converges. It has a weak singularity at $\tau = 0 $.
The line segment on the right side of the formula (\ref{formulasigmaS}) has the order of the last term:
\begin{equation}
\sum_{n=1}^{N }{\left. {{n}^{p-k}} \sin{\left( \frac{t}{{{n}^{p}}}\right) }\right.}=O(t^{(p-k)/p}).
\label{sumO}
\end{equation}
The ratio of the orders of the right-hand side in the formula (\ref{integral}) and the right-hand side of the formula (\ref{sumO}) is $t^{1/p}$. We denote $\alpha = (k-1) / p $, then the asymptotics of the series $\sigma_ {s} $ can be represented as:
$$
\sigma_{s}\sim\frac{t^{1-\alpha}}{p} 
\int_{0}^{\frac{\pi }{2}}\tau^{\alpha-2} \sin(\tau)d\tau.  
$$
Consider the sum $\sigma_ {c} $. We denote $ M = [(t / \ pi) ^ {1 / p}] $. Let's represent the sum $\sigma_{c} $ as:
\begin{equation}
\sigma_c=-\frac{1}{2}\sum_{n=1}^M n^{p-k}+\frac{1}{2}\sum_{n=1}^M n^{p-k}\cos\left(\frac{t}{n^p}\right)-\sum_{n=M+1}^\infty n^{p-k}\sin^2\left(\frac{t}{2 n^p}\right).
\label{formulaForSigmaC2}
\end{equation}
Here, as $t\to\infty$:
\begin{eqnarray*}
\frac{1}{2}\sum_{n=1}^M n^{p-k}\sim \frac{1}{2}\frac{1}{p-(k-1)}\left(\frac{t}{\pi}\right)^{1-(k-1)/p},
\\
-\frac{1}{2}\sum_{n=1}^M n^{p-k}\cos\left(\frac{1}{n^p}\right)\sim O\left(\left(\frac{t}{\pi}\right)^{1-k/p}\right).
\end{eqnarray*}

The $t$ asymptotics of the series in the formula (\ref{formulaForSigmaC2}) can be obtained in the same way as for the series from (\ref{formulasigmaS}):

$$
\sum_{n=M+1}^\infty n^{p-k}\sin^2\left(\frac{1}{2 n^p}\right)\sim \int_{M}^\infty n^{p-k}\sin^2\left(\frac{1}{2 n^p}\right)dn.
$$
replace in the integral: $n=(\frac {t}{\tau})^{1/p}$, when:
$$
\int_{M}^\infty n^{p-k}\sin^2\left(\frac{1}{2 n^p}\right)dn=\frac{{t}^{1-\frac{k-1}{p}}}{p} \int_0^{\pi}\tau^{\frac{k-1}{p}-2} \sin^2( \tau/2)d\tau 
$$
Then the $\sigma_c$ asymptotic can be represented as:
$$
\sigma_c\sim -\left(\frac{1}{2\pi^{1-\alpha}p(1-\alpha)} + 
\frac{1}{p}\int_0^{\pi}\tau^{\alpha-2} \sin^2( \tau/2)d\tau\right)t^{1-\alpha}.
$$
\begin{figure}
	\begin{subfigure}{0.45\textwidth}
		\includegraphics[scale=0.3]{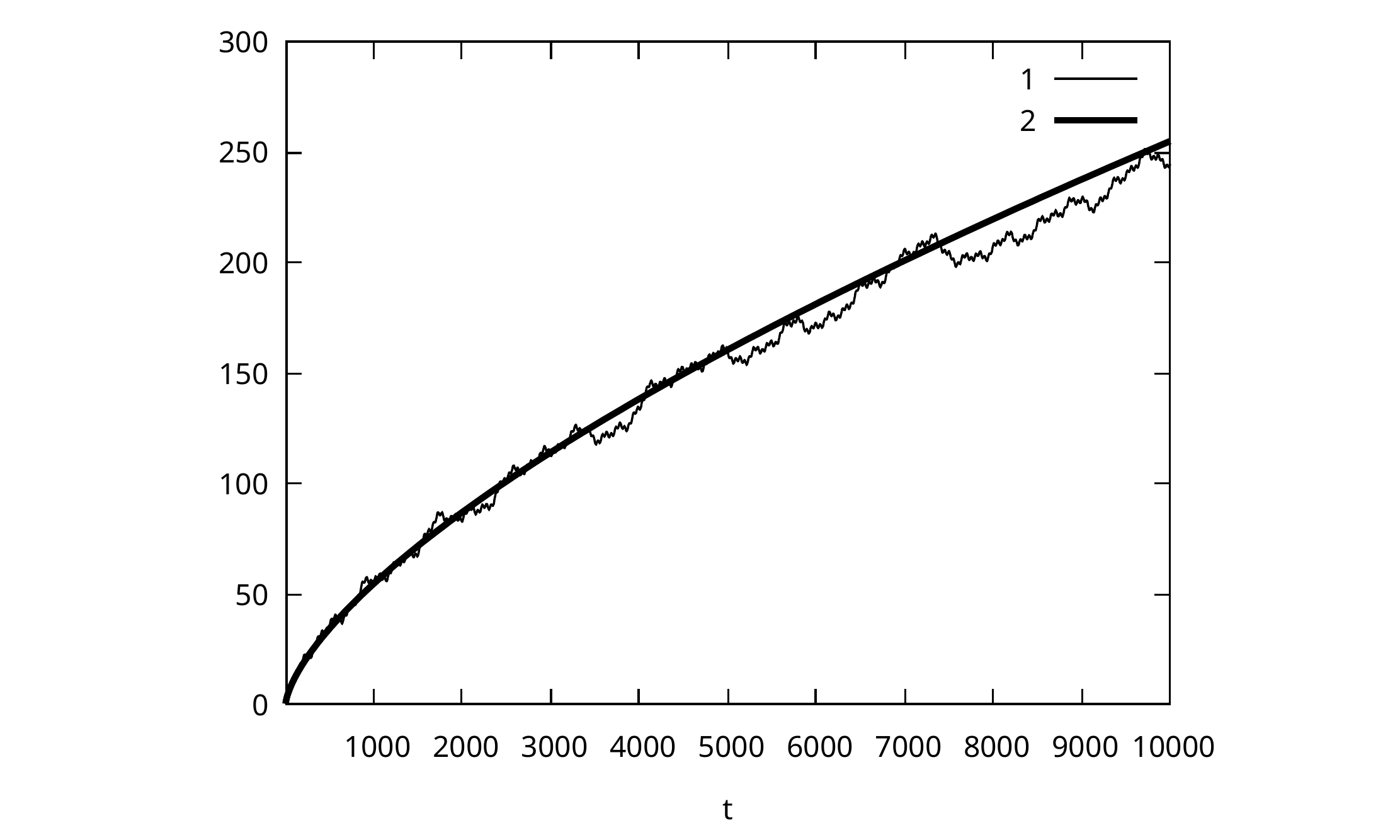}
		\caption{The coefficient for the function  $\sin(t)$ in the formula  (\ref{resh}) for the solution. Curve 1 is calculated from the first 200 terms of the series, curve 2 is calculated from the asymptotics of the main term as $t \to \infty$.}
\end{subfigure}
\hspace{0.5cm}
\begin{subfigure}{0.45\textwidth}
	\includegraphics[scale=0.3]{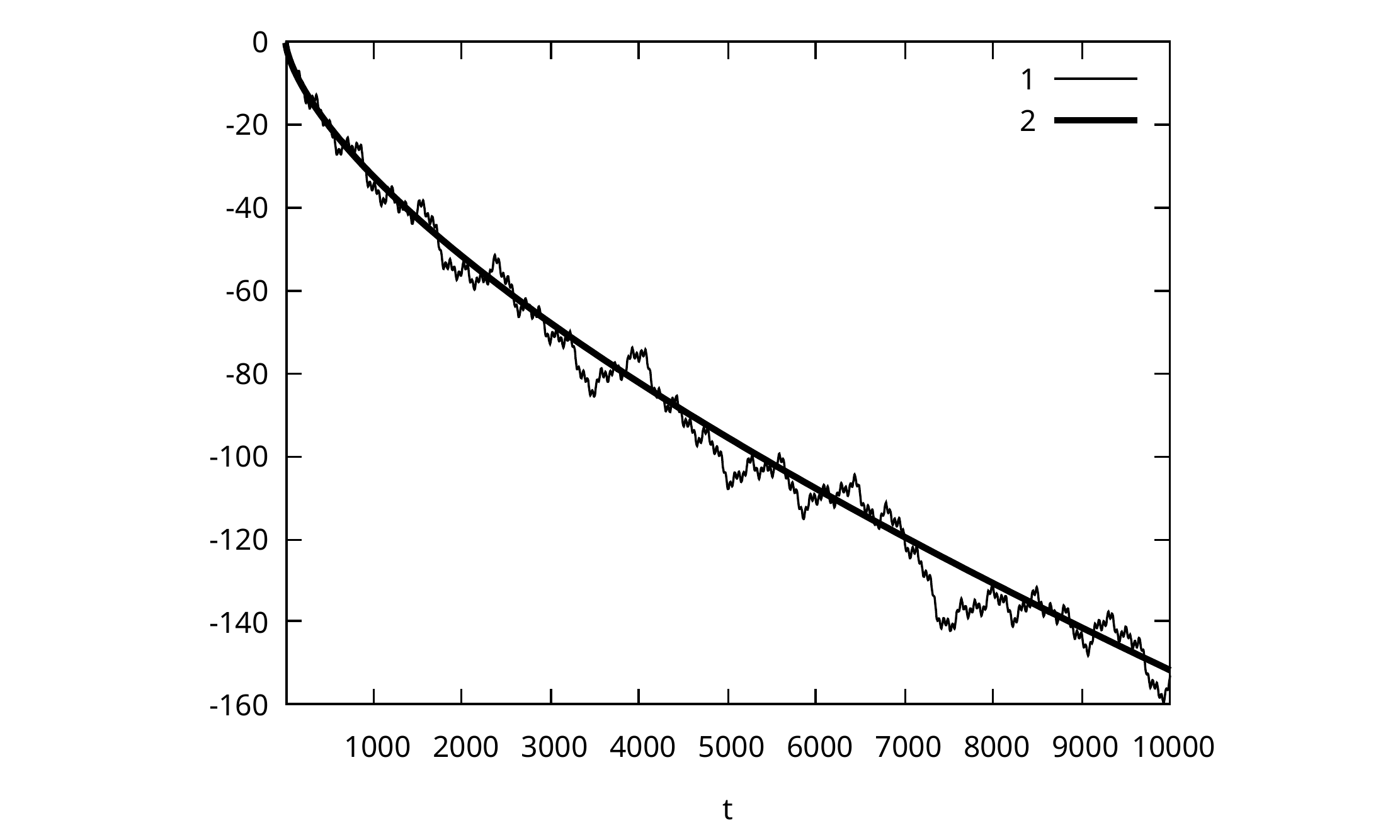}
	\caption{The coefficient for the function $\cos(t)$ in the formula  (\ref{resh}) for the solution. Curve 1 is calculated from the first 200 terms of the series, curve 2 is calculated from the asymptotics of the main term as $t \to \infty$.}
\end{subfigure}

\end{figure}

\section{The asymptotics of the subresonance solution}
\label{secAsymptoticsOfSubresonantSolution}
From the formula for the solution (\ref{resh}) and the calculations of the section  \ref{secSymptoticsOfSeries}, it follows that the asymptotic behavior of the subresonant solution of the equation  (\ref{eqLinOsc1}) has the form:
\begin{equation}
	u\sim (C_s\sin(t)+C_c\cos(t))t^{1-\alpha},
	\label{formulaForAsymptoticSolution}
\end{equation} 
where 
$$
C_s=\frac{1}{p} 
\int_{0}^{\frac{\pi }{2}}\tau^{\alpha-2} \sin(\tau)d\tau, \quad
C_c=-\left(\frac{1}{2\pi^{1-\alpha}p(1-\alpha)} + 
\frac{1}{p}\int_0^{\pi}\tau^{\alpha-2} \sin^2( \tau/2)d\tau\right).  
$$
Denote:
$$
A_{\alpha}=p\sqrt{C_s^2+C_c^2},\quad 
\phi_{\alpha}=\arctan\left(\frac{C_c}{C_s}\right),
$$
then the asymptotics of the solution  (\ref{formulaForAsymptoticSolution}) can be represented as:
\begin{equation}
	u\sim \frac{1}{p}A_\alpha t^{1-\alpha}\sin(t+\phi_{\alpha}).
	\label{formulaForAsymptotics}
\end{equation}

Note that for  $\alpha\to0$ , the growth rate over  $t$ of the amplitude of the subresonance solution tends to be linear, and   $\phi_{\alpha}\to0$. That is, the oscillating part of the solution tends to the oscillating part in the resonant solution.

The asymptotics for  $t\to\infty$ of the solution of the Cauchy problem for the equation (\ref{eqLinOsc1})  with an almost-periodic right-hand side (\ref{formulaForF}) has the form (\ref{formulaForAsymptotics}), where $\alpha=(k-1)/p$.

\end{document}